\documentclass{article}
\title{Coherence for Categorified Operadic Theories}
\author{Miles Gould\\ University of Glasgow}
\date{\today}
\usepackage{amsmath}
\usepackage{amsthm}
\DeclareSymbolFont{AMSb}{U}{msb}{m}{n}
\DeclareMathSymbol{\natural}{\mathbin}{AMSb}{"4E}
\DeclareMathSymbol{\integer}{\mathbin}{AMSb}{"5A}
\DeclareMathSymbol{\real}{\mathbin}{AMSb}{"52}
\DeclareMathSymbol{\rational}{\mathbin}{AMSb}{"51}
\DeclareMathSymbol{\I}{\mathbin}{AMSb}{"49}
\DeclareMathSymbol{\complex}{\mathbin}{AMSb}{"43}

\newcommand{\toletter}[1]{
  \mathop{\stackrel{\scriptstyle{#1}}{\longrightarrow}}
}
\newcommand{\defterm}[1]{\textbf{#1}}
\newcommand{\Cat}{\textbf{Cat}}
\newcommand{\Set}{\textbf{Set}}
\newcommand{\Operad}{\textbf{Operad}}
\newcommand{\CatOperad}{\textbf{Cat-Operad}}

\theoremstyle{plain}
\newtheorem{theorem}{Theorem}[section]
\newtheorem{lemma}[theorem]{Lemma}

\theoremstyle{definition}
\newtheorem{defn}[theorem]{Definition}

\allowdisplaybreaks
\usepackage{epsfig}
\usepackage[all]{xy}
\xyoption{2cell}
\UseTwocells
\UseHalfTwocells
\UseCompositeMaps
\CompileMatrices
\newcommand{\st}[1]{{\rm\textbf{st}}(\ensuremath{#1})}
\newcommand{\stfunc}{{\rm\textbf{st} }}
\newcommand{\seq}{_\bullet}
\newcommand{\udot}{^\bullet}
\newcommand{\dseq}{\seq\udot}
\newcommand{\adot}{a\seq}
\newcommand{\cmp}[1]{\langle#1\rangle} 
\newcommand{\ij}{_i^j}
\newcommand{\ijseq}{_{i\bullet}^j}
\newcommand{\xylabel}[1]{\ensuremath{\xymatrix{*+[o][F-]{#1}}}}

\newcommand{\WkPCat}{\ensuremath{\mbox{\rm\textbf{Wk}-}P\mbox{-\rm\textbf{Cat}}}}
\newcommand{\StrPCat}{\ensuremath{\mbox{\rm\textbf{Str}-}P\mbox{-\rm\textbf{Cat}}}}
\newcommand{\prodkn}[1]{#1_{k_1} \times \dots \times #1_{k_n}}
\newcommand{\midlabel}[1]{
	\xymatrixrowsep{4pc}
	\xymatrixcolsep{4pc}
	\xymatrix{ {} \ar @{}[d]^{#1} \\ {} }
}
\newcommand{\midequals}{\midlabel{=}}
\begin{document}

\maketitle

\begin{abstract}
It has long been known that every weak monoidal category $A$ is equivalent via
monoidal functors and monoidal natural transformations to a strict monoidal
category \st A. We generalise the definition of weak monoidal category to give a
definition of weak $P$-category for any strongly regular (operadic) theory
$P$, and show that every weak $P$-category is equivalent via $P$-functors and
$P$-transformations to a strict $P$-category. This strictification functor
is then shown to have an interesting universal property.
\end{abstract}

\section{Introduction}

Many definitions exist of categories with some kind of ``weakened'' algebraic
structure, in which the defining equations hold only up to coherent
isomorphism. The paradigmatic example is the theory of weak monoidal
categories, as presented in \cite{catwork}, but there are also definitions of
categories with weakened versions of the structure of groups \cite{baez+lauda},
Lie algebras \cite{baez+crans}, crossed monoids \cite{crossedmonoid},
sets acted on by a monoid \cite{modulecat}, rigs \cite{laplaza}, and
others. A general definition of such categories-with-structure is
obviously desirable, but hard in the general case. In this paper, we restrict
our attention to the case of strongly regular theories (equivalently, those
given by non-symmetric operads) and present possible definitions of weak
$P$-category and weak $P$-functor for any non-symmetric operad $P$. In support
of this definition, we present a generalisation of Joyal and Street's result
from \cite{j+s} that every weak monoidal category is monoidally equivalent to
a strict monoidal category.

The idea is to consider the strict models of our theory as algebras for an
operad, then to obtain the weak models as (strict) algebras for a weakened
version of that operad (which will be a \Cat-operad).
In particular, we do not make use of the pseudo-algebras of Blackwell, Kelly
and Power (for which see \cite{bkp}).
We weaken the operad
using a similar approach to that used in Penon's definition of $n$-category:
see \cite{penon}, or \cite{cheng+lauda} for a non-rigorous summary.  The weak
$P$-categories obtained are the ``unbiased'' ones: for instance, if $P$ is the
terminal operad (whose strict algebras are monoids), then the weak
$P$-categories will have tensor products of all arities, not just 0 and 2. 

In section \ref{defsect}, we present our definitions of weak $P$-category and
weak $P$-functor. In section \ref{thmsect} we extend Joyal and Street's proof
(or rather, Leinster's unbiased version) to the more
general case of weak $P$-categories. In section \ref{stsect}, we examine the
strictification functor defined in section \ref{thmsect}, and show that it has
an interesting universal property. In section \ref{furthersect}, we explain
why our approach cannot be straightforwardly extended to theories which are
not given by operads, and outline some of the approaches that could be taken
to deal with this.

\section{Weak $P$-categories}
\label{defsect}

By a \defterm{plain} operad, we mean what is elsewhere called a
``non-symmetric'' or ``non-$\Sigma$'' operad, that is one with no symmetric
group action defined on it. Unless stated otherwise, all our operads are
plain. The category of plain operads and their morphisms is \Operad.
We are also interested in plain operads enriched in \Cat:
a \defterm{\Cat-operad} is a sequence of categories $Q(0), Q(1),
\dots$, a family of composition functors $\circ:Q(n) \times Q(k_1)
\times \ldots \times Q(k_n) \to Q(\sum k_i)$ and an identity 
$1_Q \in Q(1)$, satisfying the usual associativity and unit axioms (as given,
for instance, in \cite{hohc} section 2.2).
Since operads can be thought of as one-object multicategories, we shall refer
to the objects of the categories $Q(i)$ as \defterm{1-cells} and the arrows of
these categories as \defterm{2-cells} of $Q$. 
$\Cat$-operads and their morphisms form a category: we call this
\CatOperad.
We could consider operads enriched in any symmetric monoidal category
$\mathcal{V}$, but here we are only concerned with the cases $\mathcal{V} =
\Cat$ or \Set.

Let $\mathcal{V} = \Cat$ or \Set, and $Q$ be a $\mathcal{V}$-operad.
If $A \in \mathcal{V}$, we shall write $Q \circ A$ for the
coproduct $\coprod_{n \in \natural} Q(n) \times A^n$ (this
notation was introduced by Kelly in \cite{kelly} for clubs). An
\defterm{algebra} for a $\mathcal{V}$-operad $Q$ is an object $A \in
\mathcal{V}$ and an arrow $h: Q~\circ~A \to A$ which commutes with composition
in $Q$ and such that $h(1_Q,-)$ is the identity on $A$.

Throughout, let $P$ be a plain (\Set-)operad.

A \defterm{strongly regular} algebraic theory is one that can be presented
using equations that use the same variables in the same order on both sides,
with each variable appearing only once on each side.
For instance, the theory of monoids is strongly regular, as is the theory of
sets acted on by a given monoid $M$. The theory of commutative monoids is not
strongly regular (intuitively, because of the equation $a\cdot b = b\cdot a$)
and the theory of groups is not strongly regular (again intuitively, because
of the equation $g\cdot g^{-1} = 1$). It can be shown, for instance as in
\cite{hohc} section C.1, that the strongly regular theories are exactly those
given by plain operads, in the sense that the models of a strongly regular
theory $T$ are exactly the algebras for a plain operad $P_T$.

Plain operads are algebras for a straightforward multi-sorted algebraic
theory, so (by standard arguments from universal algebra) there is an adjunction
\[
\xymatrix{\Set^\natural \ar@<1.2ex>[r]^F & \Operad \ar@<1.2ex>[l]^U_\bot}
\]
The left adjoint is given by taking labelled trees, as described in
\cite{hohc} section 3.2. Let $D: \Operad \to \CatOperad$ be the functor which
takes discrete categories aritywise; i.e., $DP(n)$ is the discrete category on
the set $P(n)$.

\begin{defn}
\label{weakeningdef}
The \defterm{weakening of $P$}, Wk($P$), is the \Cat-operad with the same
1-cells as $FUP$, and the unique 2-cell structure such that the extension of
the counit is a map of \Cat-operads and is full and faithful aritywise.
\end{defn}

\begin{figure}[h]
\centerline{
	\epsfxsize=3in
	\epsfbox{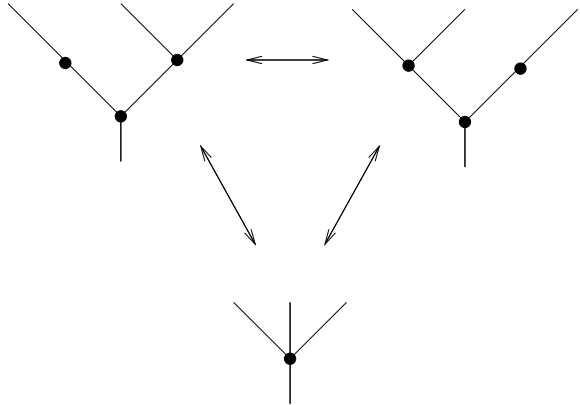}
}
\caption{Part of Wk($P$)(3) with $P = 1$}
\label{Wk(1)pic}
\end{figure}

More concretely, take $FUP$, and, for any $A, B \in FUP(n)$, place an arrow $A
\rightarrow B$ iff $\varepsilon(A) = \varepsilon(B)$ (where $\varepsilon$ is
the counit of the adjunction $F \dashv U$). The composite of two
arrows $A \rightarrow B \rightarrow C$ is the unique arrow $A \rightarrow C$.
In particular, the arrows $A \rightarrow B$ and $B \rightarrow A$ are
inverses.  See Fig. \ref{Wk(1)pic}.

This leads us immediately to the following definition:

\begin{defn}
A \defterm{weak $P$-category} is an algebra for Wk($P$).
\end{defn}

In the case $P = 1$, this reduces exactly to Leinster's definition of unbiased
monoidal category in \cite{hohc} section 3.1. There, two 1-cells $\phi$ and
$\psi$ have the same image under $\varepsilon$ iff they have the same arity, so
the categories Wk(1)$(i)$ are indiscrete.  We refer to the image under $h$ of
a map $q \to q'$ in Wk($P$) as $\delta_{q,q '}$. This is clearly a natural
transformation $h(q,-) \to h(q',-)$. As a special case, we write
$\delta_q$ for $\delta_{q, \varepsilon(q)}$.

It is not currently clear how this definition of weak $P$-category
relates to that of Hu and Kriz (as introduced in \cite{h+k}, and
elucidated in Fiore's paper \cite{fiore}), as I have only recently become
aware of their work.

\begin{defn}
A \defterm{strict $P$-category} is an algebra for $DP$, or equivalently a weak
$P$-category in which every component of $\delta$ is an identity arrow.
\end{defn}

\begin{defn}
Let $(A, h)$ and $(B, h')$ be weak $P$-categories. A \defterm{weak
$P$-functor} from $(A, h)$ to $(B, h')$ is a pair $(G, \psi)$, where $G:A
\to B$ is a functor and $\psi$ is a sequence of natural isomorphisms
$\psi_i : h'_i (1 \times G^i) \to G h_i$, such that the following diagrams
commute, for all $n, k_1, \dots, k_n \in \natural$:
\begin{equation}
\label{weakfunctordef}
\xymatrixrowsep{4pc}
\xymatrixcolsep{4pc}
\xymatrix{
	{}
		\ar[d]_{1\times 1^n \times G^{\sum k_i}}
		\ar[r]^{\prodkn{h}}
		\drtwocell\omit{^*{!(-1,-1.5)\object{\prodkn{\psi}}}}
	& {}
		\ar[d]|{1\times G^n}
		\ar[r]^{h_n}
		\drtwocell\omit{^*{!(-1,-1.5)\object{\psi_n}}}
	& {}
		\ar[d]^G
		\\
	{}
		\ar[r]_{\prodkn{h'}}
	& {}
		\ar[r]_{h'}
	& {} \\
}
\midequals
\xymatrixrowsep{4pc}
\xymatrixcolsep{4pc}
\xymatrix{
	{}
		\ar[d]_{1\times 1^n \times G^{\sum k_i}}
		\ar[r]^{h_{\sum k_i}}
		\drtwocell\omit{^*{!(-1,-1.5)\object{\psi_{\sum k_i}}}}
	& {}
		\ar[d]^G
		\\
	{}
		\ar[r]_{h'_{\sum k_i}}
	& {} \\
}
\end{equation}
\begin{equation}
\label{weakfunctordef2}
\xymatrixrowsep{3pc}
\xymatrixcolsep{3pc}
\xymatrix{
	G a \ar[r]^{\delta'_1} \ar[d]_1
	& h'(1_P, G a) \ar[d]^{\psi_1} \\
	G a \ar[r]_{G \delta_1}
	& G h(1_P, a)
}
\end{equation}
This definition is a natural generalisation of the definition of weak unbiased
monoidal functor given in \cite{hohc} section 3.1.
\end{defn}
\begin{defn}
Let $(F, \phi)$ and $(G, \psi)$ be weak $P$-functors $(A, h) \to (B, h')$. A
\defterm{$P$-transformation} $\sigma: (F, \phi) \to (G, \psi)$ is a natural
transformation 
\[
\xymatrix{ A \rtwocell^F_G{\sigma} & B }
\]
\end{defn}
such that
\begin{equation}
\xymatrixrowsep{4pc}
\xymatrixcolsep{4pc}
\xymatrix{
	Wk(P) \circ A \ar[r]^h
	\dtwocell^{<2>{1 \circ G}}_{<2>{1 \circ F}} {^\sigma} &
	A \ar[d]^G \\
	Wk(P) \circ B \ar[r]^{h'} &
	B
	\ultwocell\omit{\psi}
}
\midequals
\xymatrix{
	Wk(P) \circ A \ar[r]^h \ar[d]_{1 \circ F} &
	A \dtwocell^G_F {^\sigma} \\
	Wk(P) \circ B \ar[r]^{h'} &
	B
	\ultwocell\omit{\phi}
}
\end{equation}
Note that there is only one possible level of strictness here. There is a
2-category, \WkPCat, whose objects are weak $P$-categories, whose 1-cells are
weak $P$-functors, and whose 2-cells are $P$-transformations. Similarly, there
is a 2-category \StrPCat\ of strict $P$-categories, strict $P$-functors, and
$P$-transformations, which can be considered a sub-2-category of \WkPCat.

\begin{lemma}
\label{inv <=> P-inv}
A $P$-transformation $\sigma: (F, \phi) \to (G, \psi)$ is invertible as a
$P$-transformation if and only if it is invertible as a natural
transformation.
\end{lemma}
\begin{proof}
``Only if'' is obvious: we concentrate on ``if''. It's enough to show that
$\sigma^{-1}$ is a $P$-transformation, which is to say that
\begin{equation}
\xymatrix{
	h(q, G\adot)
		\ar[r]^\psi
		\ar[d]_{h(q,\sigma^{-1}_{\adot})} &
	Gh(q,\adot)
		\ar[d]^{\sigma^{-1}_{h(q,\adot)}} \\
	h(q, F\adot)
		\ar[r]^\phi &
	Fh(q,\adot)
}
\end{equation}
commutes for all $(q,\adot) \in Wk(P) \circ A$, and this follows
from the fact that ${\sigma_{h(q,\adot)}} \circ \phi = \psi \circ
{h(q,\sigma_{\adot})}$.

\end{proof}

\section{Main Theorem}
\label{thmsect}

Let $P$ be a plain operad, and $Q = \mbox{Wk}(P)$, with $\pi:Q \to P$ the
projection map. We write composition in $P$ as $p \circ (p_1,\dots,p_n)$, and
composition in $Q$ as $q\cmp{q_1,\dots,q_n}$. We also adopt the ${}\seq$
notation from chain complexes and write, for instance, $p\seq$ for a sequence
of objects in $P$ and $p\seq\udot$ for a double sequence. Let $Q \circ A
\toletter{h} A$ be a weak $P$-category. We construct a strict $P$-category \st
A and a weak $P$-functor $(F,\phi):\st A \to A$, and show that it is an
equivalence of weak $P$-categories.

In fact, \stfunc is functorial, and is left adjoint to the forgetful
functor $\StrPCat\linebreak[0] \to \WkPCat$ (see \S \ref{stsect}).
The theorem then says that the
unit of this adjunction is pseudo-invertible, and that the strict
$P$-categories and strict $P$-functors form a weakly coreflective
sub-2-category of \WkPCat.

\begin{defn}
Let $P$, $Q$, $A$ be as above. The \defterm{strictification of $A$}, \st
A, is defined as follows:
\begin{itemize}
\item An object of \st A is an object of $\coprod P(i) \times A^i$.
\item An arrow $(p, \adot) \to (p', \adot')$ in \st A
is an arrow $h(p, \adot) \to h(p', \adot')$ in $A$.
Composition and identities are as in $A$.
\end{itemize}

We define an action $h'$ of $Q$ on \st A as follows:
\begin{itemize}
\item On objects, $h'$ acts by $h'(q, (p, \adot)\udot)
= (\pi(q\cmp {p\udot}), \adot\udot)$.
\item Let $f_i :(p_i, a_i) \to (p_i', a_i')$ for $i = 0, \dots,
  n$. Then $h'(p, f\seq)$ is the composite
\[
\begin{array}{rcl}
h(p\circ(p\seq), \adot)& \toletter{\delta_{p\cmp {p\seq}}^{-1}}&
        h(p\cmp {p\seq}, \adot) = h(p, h(p_0, a_0), \dots , h(p_n, a_n))\\
    &\toletter{h(p, f\seq)}
        &h(p, h(p_0', a_0'), \dots, h(p_n', a_n'))
        = h(p\cmp {p\seq'}, \adot')\\
    &\toletter{\delta_{p\cmp {p\seq'}}}
        &h(p\circ(p\seq'), \adot').
\end{array}
\]
\end{itemize}
\end{defn}

\begin{lemma}\st A is a strict $P$-category.\end{lemma}

\begin{proof}
The identity and strictness conditions are obvious, as is the associativity
of the action on objects.  We must show that the action on arrows is
associative. Let $f_i^j:(p\ij, a\ijseq) \to (q\ij, b\ijseq)$, $\sigma
\in Q(n)$, and $\tau_i \in Q(k_i)$ for $j = 1, \dots, k_i$ and $i = 1,\dots,n$.
We wish to show that $h'(\sigma\circ(\tau\seq), f\dseq) = h'(\sigma,
h'(\tau_1, f_1\udot), \dots, h'(\tau_n, f_n\udot))$.

The LHS is
\[
\begin{array}{rcl}
h(\sigma\circ(\tau\seq)\circ(p\dseq), \adot\udot)&
	\toletter{\delta_{\sigma\circ(\tau_i)\cmp{p\dseq}}^{-1}}&
h(\sigma\circ(\tau\seq), h(p_1^1, a_{1\bullet}^1), \dots, h(p_n^{k_n},
a_{n\bullet}^{k_n}))\\
	&\toletter{h(\sigma\circ(\tau\seq), f\dseq)}
&h(\sigma\circ(\tau\seq), h(q_1^1, b_{1\bullet}^1), \dots, h(q_n^{k'_n},
b_{n\bullet}^{k'_n}))\\
	&\toletter{\delta_{\sigma\circ(\tau\seq)\cmp{q\dseq}}}
&h(\sigma\circ(\tau\seq)\circ(q\dseq), b\dseq).
\end{array}
\]

The RHS is 
\[
\begin{array}{rcl}
h(\sigma\circ(\tau\seq)\circ(p\dseq), \adot\udot)&
	\toletter{\delta_{\sigma\cmp{\tau_i\circ(p\dseq)}}^{-1}}&
h(\sigma, h(\tau_1\circ(p_1\udot), a_{1\bullet}\udot), \dots,
		h(\tau_n\circ(p_n\udot), a_{n\bullet}\udot)) \\
	&\toletter{h(\sigma, h'(\tau\seq, f\dseq))}
&h(\sigma, h(\tau_1\circ(q_1\udot), b_{1\bullet}\udot), \dots,
		h(\tau_n\circ(q_n\udot), b_{n\bullet}\udot))\\
	&\toletter{\delta_{\sigma\cmp{\tau_i\circ(p\dseq)}}}
&h(\sigma\circ(\tau\seq)\circ(q\dseq), b\dseq),
\end{array}
\]
where each $h'(\tau_i, f_i\udot)$ is
\[
\begin{array}{rcl}
h(\tau_i\circ(p_i\udot),a_i\udot))
	&\toletter{\delta_{\tau_i\cmp{p_i\udot}}^{-1}}
& h(\tau_i, h(p_i^1, a_{i\bullet}^1), \dots, h(p_i^{k_i}, a_{i\bullet}^{k_i}))\\
	&\toletter{h(\tau_i, f_i\udot)}
& h(\tau_i, h(q_i^1, b_{i\bullet}^1), \dots, h(q_i^{k_i}, b_{i\bullet}^{k_i}))\\
	&\toletter{\delta_{\tau_i\cmp{p_i\udot}}}
&h(\tau_i\circ(q_i\udot),b_i\udot).
\end{array}
\]

So the equation holds if the following diagram commutes:

\[
\xymatrix @+10pt {
& h(\sigma\circ(\tau\seq)\circ(p\dseq), \adot\udot)
	\ar[dl]_{\delta_{\sigma\circ(\tau_i)\cmp{p\dseq}}^{-1}}
	\ar[dr]^{\delta_{\sigma\cmp{\tau_i\circ(p\dseq)}}^{-1}}
	\ar[d]|{\delta_{\sigma\cmp{\tau\seq}\cmp{p\dseq}}^{-1}} \\
h(\sigma\circ(\tau\seq), h(p\dseq, a\dseq))
	\ar[d]|{h(\sigma\circ(\tau\seq), f\dseq)}
& h(\sigma, h(\tau\seq\circ(p\dseq),a\dseq))
	\ar[l]^{\delta_{\sigma\cmp{\tau\seq}}}
	\ar[d]|{h(\sigma, h(\tau_i, f_i\udot))}
	\ar @{} [dl]|*+[o][F-]{1}
	\ar @{} [dr]|*+[o][F-]{2}
& h(\sigma\circ(\tau\seq), h(p\dseq, a\dseq))
	\ar[l]^{h(\sigma, \delta_{\tau\seq\cmp{p\dseq}}^{-1})}
	\ar[d]|{h(\sigma, h'(\tau\seq, f\dseq))}\\
h(\sigma\circ(\tau\seq), h(q\dseq, b\dseq))
	\ar[dr]_{\delta_{\sigma\circ(\tau\seq)\cmp{q\dseq}}}
& h(\sigma, h(\tau\seq\circ(q\dseq), b\dseq))
	\ar[d]|{\delta_{\sigma\cmp{\tau\seq\circ(q\dseq)}}}
	\ar[l]^{\delta_{\sigma\cmp{\tau\seq}}}
	\ar[r]_{h(\sigma, \delta_{\tau_i\cmp{p_i\udot}})}
& h(\sigma, h(\tau\seq\circ(q\dseq), b\dseq))
	\ar[dl]^{\delta_{\sigma\cmp{\tau\seq\circ(q\dseq)}}}\\
&h(\sigma\circ(\tau\seq)\circ(q\dseq), b\dseq)\\
}
\]
The triangles all commute because all $\delta$s are images of arrows in $Q$,
and there is at most one 2-cell between any two 1-cells in $Q$. \xylabel 2
commutes by the definition of $h'(\tau_i, f_i\udot)$, and \xylabel 1
commutes by naturality of $\delta$.
\end{proof}

\begin{lemma}
\label{equiv=>Q-equiv}
Let $Q\circ A \toletter{h} A$ and $Q \circ B \toletter{h'} B$
be weak $P$-categories, $(F, \pi):A \to B$ be a weak $P$-functor, and $(F, G,
\eta, \varepsilon)$ be an adjoint equivalence. Then $G$ naturally carries the
structure of a weak $P$-functor, and $(F, G, \eta, \varepsilon)$ is an adjoint
equivalence in \WkPCat.
\end{lemma}

\begin{proof}
We want a sequence $(\psi\seq)$ of natural transformations:
\[
\xymatrixrowsep{3pc}
\xymatrixcolsep{3pc}
\xymatrix{
	Q(i)\times B^i
		\ar[d]_{1\times G^i}
		\ar[r]^{h'_i}
		\drtwocell\omit{^*{!(-1,-1.5)\object{\psi_i}}}
		& B \ar[d]^G \\
	Q(i)\times A^i \ar[r]_{h_i}
		& A \\
}
\]
Let $\psi_i$ be given by $Gh'(1 \times \varepsilon^i) \circ G\pi_i^{-1} \circ
\eta_{h'}$, i.e.\
\def\epsiloncell#1{%
	\ultwocell\omit{<-0.7>*{!(-2,-1.2){1\times \varepsilon^{#1}}}}
}
\def\epsilonsumki{\epsiloncell{\sum k_i}}
\def\epsilonn{\epsiloncell{n}}
\[
\xymatrixrowsep{3.5pc}
\xymatrixcolsep{3pc}
\xymatrix{
	Q(i) \times B^i
		\ar[d]_{1\times G^i}
		\ar[r]^{h'_i}
		\drtwocell\omit{^*{!(-1,-1.5)\object{\psi_i}}}
	& B
		\ar[d]^G \\
	Q(i) \times A^i
		\ar[r]_{h_i}
	& A \\
}
\midequals
\xymatrixrowsep{1.5pc}
\xymatrixcolsep{1.2pc}
\xymatrix{
	Q(i) \times B^i
		\ar[rr]^1
		\ar[dd]_{1\times G^i}
	&& Q(i) \times B^i
		\ar[rr]^{h'_i}
	&& B
		\ar[dd]^G
		\\
	& {}
	\ultwocell\omit{<-0.7>*{!(-2,-1.8){1\times \varepsilon^{i}}}}
	&& {}
		\lltwocell\omit{*{!(-5,-4)\object{{\pi^{-1}_i}}}}
		\\
	Q(i) \times A^i
		\ar[urur]|{1 \times F^i}
		\ar[rr]_{h_i}
	&& A
		\ar[urur]^F
		\ar[rr]_1
	&& A
		\ultwocell\omit{{\eta}} \\
}
\]
We must check that $\psi$ satisfies (\ref{weakfunctordef}) and
(\ref{weakfunctordef2}). For (\ref{weakfunctordef}):
\begin{eqnarray*}
\midlabel{\mbox{LHS}} & \midequals
& \xymatrixrowsep{5pc}
\xymatrixcolsep{5pc}
\xymatrix{
	{}
		\ar[d]_{1\times G^{\sum k_i}}
		\ar[r]^{\prodkn{h'}}
		\drtwocell\omit{^*{!(-1,-1.5)\object{\prodkn{\psi}}}}
	& {}
		\ar[d]|{1\times G^n}
		\ar[r]^{h'_n}
		\drtwocell\omit{^*{!(-1,-1.5)\object{\psi_n}}}
	& {}
		\ar[d]^G
		\\
	{}
		\ar[r]_{\prodkn{h}}
	& {}
		\ar[r]_h
	& {} \\
} \\
& \midequals
& \xymatrixrowsep{2.3pc}
\xymatrixcolsep{2.3pc}
\xymatrix{
	{}
		\ar[rr]^1
		\ar[dd]_{1\times G^{\sum k_i}}
	&& {}
		\ar[rr]^{\prodkn{h'}}
	& {}
	& {}
		\ar[rr]^1
		\ar[dd]|{1\times G^n}
	&& {}
		\ar[rr]^{h'_n}
	& {}
	& {}
		\ar[dd]^G
		\\
	& {}
		\epsilonsumki
	&& {}
		\lltwocell\omit{*{!(-5,-4)\object{\prodkn{\pi^{-1}}}}}
	&& {}
		\epsilonn
	&& {}
		\lltwocell\omit{*{!(-5,-4)\object{\pi_n^{-1}}}}
	&& \\
	{}
		\ar[uurr]|{1 \times F^{\sum k_i}}
		\ar[rr]_{\prodkn{h}}
	&& {}
		\ar[uurr]|{1 \times F^n}
		\ar[rr]_1
	&& {}
		\ar[uurr]|{1 \times F^n}
		\ar[rr]_h
		\ultwocell\omit{*{!(1,0)\object{1 \times \eta^n}}}
	&& {}
		\ar[uurr]^F
		\ar[rr]_1
	&& {}
		\ultwocell\omit{{\eta}} \\
} \\
& \midequals
& \xymatrixrowsep{2.3pc}
\xymatrixcolsep{2.3pc}
\xymatrix{
	{}
		\ar[rr]^1
		\ar[dd]_{1\times G^{\sum k_i}}
	&& {}
		\ar[rr]^{\prodkn{h'}}
	& {}
	& {}
		\ar[rr]^1
	&& {}
		\ar[rr]^{h'_n}
	& {}
	& {}
		\ar[dd]^G
		\\
	& {}
		\epsilonsumki
	&& {}
		\lltwocell\omit{*{!(-5,-4)\object{\prodkn{\pi^{-1}}}}}
	&&&& {}
		\lltwocell\omit{*{!(-5,-4)\object{\pi_n^{-1}}}}
	& \\
	{}
		\ar[urur]|{1 \times F^{\sum k_i}}
		\ar[rr]_{\prodkn{h}}
	&& {}
		\ar[urur]|{1 \times F^n}
		\ar[rr]_1
	&& {}
		\ar[urur]|{1 \times F^n}
		\ar[rr]_h
		\uutwocell\omit{=}
	&& {}
		\ar[urur]^F
		\ar[rr]_1
	&& {}
		\ultwocell\omit{{\eta}} \\
} \\
& \midequals
& \xymatrixrowsep{2.3pc}
\xymatrixcolsep{2.3pc}
\xymatrix{
	{}
		\ar[rr]^1
		\ar[dd]_{1\times G^{\sum k_i}}
	&& {}
		\ar[rr]^{\prodkn{h'}}
	& {}
	& {}
		\ar[rr]^{h'_n}
	& {}
	& {}
		\ar[dd]^G
		\\
	& {}
		\epsilonsumki
	&& {}
		\lltwocell\omit{*{!(-5,-4)\object{\prodkn{\pi^{-1}}}}}
	&& {}
		\lltwocell\omit{*{!(-5,-4)\object{\pi_n^{-1}}}}
	&& \\
	{}
		\ar[urur]|{1 \times F^{\sum k_i}}
		\ar[rr]_{\prodkn{h}}
	&& {}
		\ar[urur]|{1 \times F^n}
		\ar[rr]_h
	&& {}
		\ar[urur]^F
		\ar[rr]_1
	&& {}
		\ultwocell\omit{{\eta}} \\
} \\
& \midequals
& \xymatrixrowsep{2.3pc}
\xymatrixcolsep{2.3pc}
\xymatrix{
	{}
		\ar[rr]^1
		\ar[dd]_{1\times G^{\sum k_i}}
	&& {}
		\ar[rr]^{h'_{\sum k_i}}
	& {}
	& {}
		\ar[dd]^G
		\\
	& {}
		\epsilonsumki
	&& {}
		\lltwocell\omit{*{!(-5,-4)\object{\pi_{\sum k_i}^{-1}}}}
	& \\
	{}
		\ar[urur]|{1 \times F^{\sum k_i}}
		\ar[rr]_{h_{\sum k_i}}
	&& {}
		\ar[urur]^<<<<<<<<<<<<F
		\ar[rr]_1
	&& {}
		\ultwocell\omit{{\eta}} \\
} \\
& \midequals
& \xymatrixrowsep{5pc}
\xymatrixcolsep{5pc}
\xymatrix{
	{}
		\ar[d]_{1\times G^{\sum k_i}}
		\ar[r]^{h'_{\sum k_i}}
		\drtwocell\omit{^*{!(-1,-1.5)\object{\psi_{\sum k_i}}}}
	& {}
		\ar[d]^G
		\\
	{}
		\ar[r]_{h_{\sum k_i}}
	& {} \\
} \\
& = & \mbox{RHS.}
\end{eqnarray*}

For (\ref{weakfunctordef2}), consider the following diagram:
\[
\xymatrix{
	&&G b \ar[r]^{\delta_{1_Q}} \ar@/l3cm/[ddd]_{1} \ar[d]_{\eta G}
	  \ar @{} [dr]|*+[o][F-]{2} 
	  & h(1_P, G b) \ar@/r3cm/[ddd]^{\psi_1} \ar[d]^{\eta} \\
	\ar @{} [drr]|*+[o][F-]{1} 
	&&GFG b \ar[r]^{GF\delta_{1_Q}} \ar[d]_{1}
	  \ar @{} [dr]|*+[o][F-]{3} 
	  & GFh(1_P, Gb) \ar[d]^{\pi^{-1}_1}
	  \ar @{} [drr]|*+[o][F-]{5} 
	  \\
	&&GFG b \ar[r]_{G\delta'_{1_Q}} \ar[d]_{G\varepsilon}
	  \ar @{} [dr]|*+[o][F-]{4} 
	  & G h'(1_P, FGb) \ar[d]^{Gh'(1_P, \varepsilon)} &&\\
	&&Gb \ar[r]_{G\delta'_{1_Q}}
	  & Gh'(1_P, b)
}
\]
(\ref{weakfunctordef2}) is the outside of the diagram.
\xylabel{1} commutes by the triangle identities.
\xylabel{2} commutes by naturality of $\eta$.
\xylabel{3} commutes since $(F, \pi)$ is a $P$-functor.
\xylabel{4} commutes by naturality of $\delta$.
\xylabel{5} is the definition of $\psi$.
Hence the whole diagram commutes, and $(G, \psi)$ is a $P$-functor.

To see that $(F, G, \eta, \varepsilon)$ is a $P$-equivalence, it is now enough to
show that $\eta$ and $\varepsilon$ are $P$-transformations, since they satisfy the
triangle identities by hypothesis.

Write $(GF, \chi) = (G, \psi) \circ (F, \pi)$.
We wish to show that $\eta$ is a $P$-transformation $(1,1) \to (GF, \chi)$.
Each $\chi_{q, \adot}$ is the composite
\[
\xymatrix{
	h(q, GF\adot) \ar[r]^\psi
	& G h(q, F\adot) \ar[r]^{G\pi}
	& GFh(q, \adot)
}
\]
Applying the definition of $\psi$, this is
\[
\centerline{
\xymatrix{
	h(q, GF\adot) \ar[r]^\eta
	& GF h(q, GF\adot) \ar[r]^{G\pi^{-1}}
	& G h(q, FGF\adot) \ar[r]^{Gh_q \varepsilon F}
	& Gh(q, F \adot)   \ar[r]^{G\pi}
	& GFh(q, \adot)
}
}
\]

The axiom on $\eta$ is the outside of the diagram
\[
\centerline{
\xymatrix{
	h(q, \adot) \ar[rrrr]^1 \ar[dd]_{h(q, \eta)} \ar[dr]^\eta
	&&&& h(q, \adot) \ar[dd]^\eta \\
	& GFh(q, \adot) \ar[r]^{G\pi^{-1}} \ar[d]|{GFh(q, \eta)}
	\ar @{} [l]|*+[o][F-]{1} 
	\ar @{} [dr]|*+[o][F-]{2} 
	& Gh(q, F\adot) \ar[d]|{Gh(q, F\eta)} \ar[dr]^1
	\ar @{} [rr]|*+[o][F-]{3} & & \\ 
	h(q, GF\adot) \ar[r]^\eta
	& GFh(q, GF\adot) \ar[r]^{G\pi^{-1}}
	& Gh(q, FGF\adot) \ar[r]^{Gh(q, \varepsilon F)}
	& Gh(q, F\adot) \ar[r]^{G\pi}
	& GFh(q, \adot)
}
}
\]
\xylabel{1} commutes by naturality of $\eta$, \xylabel{2} commutes by
naturality of $\pi^{-1}$, and \xylabel{3} commutes since $G\pi \circ G\pi^{-1}
= G(\pi \circ \pi^{-1}) = G1 = 1G$. The triangle commutes by the triangle
identities. So the whole diagram commutes, and $\eta$ is a $P$-transformation.
By Lemma \ref{inv <=> P-inv}, $\eta^{-1}$ is also a $P$-transformation.
Similarly, $\varepsilon$ and $\varepsilon^{-1}$ are $P$-transformations.

\end{proof}

\begin{theorem}
\label{mainthm}
Let $Q\circ A \toletter{h} A$ be a weak $P$-category. Then $A$
is equivalent to \st A via weak $P$-functors and $P$-transformations.
\end{theorem}

\begin{proof}
Let $F: \st A \to A$ be given by $F(p, \adot) = h(p, \adot)$ and
identification of maps. This is certainly full and faithful, and it's
essentially surjective on objects because $\delta^{-1}_{1_Q}: h(1_P, a)
\to a$ is an isomorphism. It remains to show that $F$ is a weak $P$-functor.

We must find a sequence $(\phi_i : h_i (1 \times F^i) \to Fh')$ of natural
transformations satisfying equations (\ref{weakfunctordef}) and
(\ref{weakfunctordef2}).
We can take $(\phi_i)_{q, (p\seq, a\dseq)} = (\delta_{q \cmp{p\seq}})_{a\dseq}$.
For (\ref{weakfunctordef}), we must show that

\[
\xymatrixrowsep{4pc}
\xymatrixcolsep{4pc}
\xymatrix{
	{}
		\ar[d]_{1\times F^{\sum k_i}}
		\ar[r]^{\prodkn{h'}}
		\drtwocell\omit{^*{!(-1,-1.5)\object{\prodkn{\phi}}}}
	& {}
		\ar[d]|{1\times F^n}
		\ar[r]^{h'_n}
		\drtwocell\omit{^*{!(-1,-1.5)\object{\phi_n}}}
	& {}
		\ar[d]^F
		\\
	{}
		\ar[r]_{\prodkn{h}}
	& {}
		\ar[r]_h
	& {} \\
}
\midequals
\xymatrixrowsep{4pc}
\xymatrixcolsep{4pc}
\xymatrix{
	{}
		\ar[d]_{1\times F^{\sum k_i}}
		\ar[r]^{h'_{\sum k_i}}
		\drtwocell\omit{^*{!(-1,-1.5)\object{\phi_{\sum k_i}}}}
	& {}
		\ar[d]^F
		\\
	{}
		\ar[r]_{h_{\sum k_i}}
	& {} \\
}
\]

All 2-cells in this equation are instances of $\delta$. Since there is at most
one 2-cell between two 1-cells in $Q$, the equation holds.

For (\ref{weakfunctordef2}) to hold, we must have 
\begin{equation}
\label{Fwkfctr}
\xymatrix{
	F(p,\adot) \ar[d]_1 \ar[r]^{\delta_{1_Q}}
	  & h(1_P, F(p,\adot)) \ar[d]^{\phi_{1_P}} \\
	F(p,\adot) \ar[r]_{F\delta'_{1_Q}}
	  & Fh'(1_P,(p, \adot))
}
\end{equation}
Since \st A is a strict monoidal category, $\delta' = 1$.
Apply this observation, and the definitions of $F$, $\phi$ and $h'$; then
(\ref{Fwkfctr}) becomes
\[
\xymatrix{
	h(p,\adot) \ar[d]_1 \ar[r]^{\delta_{1_Q}}
	  & h(1_P, h(p,\adot)) \ar[d]^{\delta_{1_P\cmp{p}}} \\
	h(p,\adot) \ar[r]_1
	  & h(p,\adot)
}
\]
Since there is at most one arrow between two 1-cells in $Q$, this diagram
commutes.
So $(F,\phi)$ is a weak $P$-functor.  Hence, by Lemma \ref{equiv=>Q-equiv},
$A$ is equivalent to \st A via weak $P$-functors and $P$-transformations.
\end{proof}

\section{Significance of \stfunc}
\label{stsect}
\begin{theorem}
Let $U'$ be the forgetful functor $\StrPCat \to \WkPCat$ (considering both of
these as 1-categories). Then \stfunc is left adjoint to $U'$.
\end{theorem}
\begin{proof}
For each $A \in \WkPCat$, we construct an
initial object $A \toletter{(F', \psi)} \st A$ of the comma category
$(A\downarrow U')$, thus showing that \stfunc is functorial and that $\stfunc
\dashv U'$ (and that $(F', \psi)$ is the component of the unit at $A$).  Let
$(B, h'')$ be a strict $P$-category, and
$(G, \gamma): A \to U'B$ be a weak $P$-functor. We must show that there is a
unique strict $P$-functor $H$ making the following diagram commute:
\begin{equation}
\label{A->U}
\xymatrix{
	& A \ar[ddl]_{(F', \psi)} \ar[ddr]^{(G, \gamma)} \\ \\
	U'\,\st A \ar@{-->}[rr]^{(H, {\rm id})} & & U'B
}
\end{equation}
$(F', \psi)$ is given as follows:
\begin{itemize}
\item If $a \in A$, then ${F'}(a) = (1, a)$.
\item If $f: a \to a'$ in $A$ then ${F'}f$ is the lifting of $h(1,f)$ with
source $(1, a)$ and target $(1, a')$.
\item $\psi_{(p, \adot)}$ is the lifting of $(\delta_{1_Q})_{h(p, \adot)}: h(p,
\adot) \to h(1, h(p, \adot))$ to a morphism
$h'(p, F'(a)\seq) = (p, \adot) \to (1, h(p, \adot)) = F'(h(p, \adot))$.
\end{itemize}
For commutativity of (\ref{A->U}), we must have $H(1, a) = G(a)$, and
for strictness of $H$, we must have $H(p, \adot) = h''(p, H(1, a)\seq)$. These
two conditions completely define $H$ on objects.

Now, take a morphism $f:(p, \adot) \to (p', \adot')$, which is a lifting of a
morphism $g: h(p,\adot) \to h(p',\adot')$ in $A$. $Hf$ is a morphism $h''(p,
G\adot) \to h''(p', G\adot')$: the obvious thing for it to be is the composite
\[
\xymatrix{
	h''(p, G\adot) \ar[r]^\gamma 
	& G h''(p, \adot) \ar[r]^{Gg}
	& G h''(p', \adot') \ar[r]^{\gamma^{-1}}
	& h''(p', G\adot')
}
\]
and we shall show that this is in fact the only possibility. Consider the
composite
\[
\xymatrix{
	(1, h(p, \adot)) \ar[r]^{\psi^{-1}}
	&(p, \adot) \ar[r]^{f}
	& (p', \adot) \ar[r]^{\psi}
	& (1, h(p', \adot'))
}
\]
in \st A. Composition in \st A is given by composition in $A$, so this is
equal to the lifting of $\delta_{1_Q} \circ g \circ \delta_{1_Q}^{-1} = h(1,
g)$ to a morphism $(1, h(p, \adot)) \to (1, h(p',
\adot'))$, namely $F'g$. So $f = \psi^{-1} \circ F'g \circ \psi$, and
$Hf = H\psi^{-1} \circ HF'g \circ H\psi$. By commutativity of
(\ref{A->U}), $HF' = G$ and $H\psi = \gamma$, so $Hf = \gamma^{-1} \circ
Gg \circ \gamma$ as required.

This completely defines $H$. So we have constructed a unique
$H$ which makes (\ref{A->U}) commute and which is strict. Hence $(F',
\psi): A \to U'\,\st A$ is initial in $(A\downarrow U')$, and so $\stfunc
\dashv U'$.
\end{proof}

The $P$-functor $(F, \phi):\st A \to A$ constructed in Theorem \ref{mainthm} is
pseudo-inverse to $(F', \psi)$, which we have just shown to be the
$A$-component of the unit of the adjunction \stfunc $\dashv U'$. We can
therefore say that \StrPCat\ is a weakly coreflective sub-2-category of
\WkPCat. Note that the counit is \emph{not} pseudo-invertible, so this is not
a 2-equivalence.

\section{Further Work}
\label{furthersect}
Very few interesting theories are strongly regular, so this definition is
unsatisfactory as it stands. It can be straightforwardly extended to theories
given by symmetric operads,
but to deal with the interesting cases of groups, rings, Lie algebras, etc, we
must either abandon operads and move to a more expressive formalism (for
instance that of Lawvere theories), or extend the notion of an operad until it
is sufficiently expressive. I have taken the latter approach: by allowing any
function of finite sets, and not just permutations, to act on the sets $P(i)$,
we obtain a notion of operad that is equivalent in power to clones or Lawvere
theories (as was proved by Tronin in \cite{tronin}).

However, na\"ively extending definition \ref{weakeningdef} to these
more general operads doesn't work, as weakening the theory of commutative
monoids gives the theory of strictly symmetric weak monoidal categories,
rather than that of symmetric weak monoidal categories as desired. I have been
working off and on on various other approaches, mainly concerned with
constructing Wk($P$) using some universal property, and have obtained some
interesting early results. 
\bibliographystyle{hplain}
\bibliography{strictify}

\begin{thebibliography}{10}

\bibitem{crossedmonoid}
Pierre Ageron.
\newblock Les cat\'egories mono\"idales crois\'ees.
\newblock 78th PSSL, Strasbourg, Feb 2002.
\newblock Available at http:/\slbrk www.math.unicaen.fr\slbrk
  {$\sim$}ageron\slbrk Strasbourg.dvi.

\bibitem{baez+crans}
John~C. Baez and Alissa~S. Crans.
\newblock Higher-dimensional algebra {VI}: Lie 2-algebras.
\newblock {\em Theory and Applications of Categories}, 12:492--538, 2004,
  math.QA/0307263.

\bibitem{baez+lauda}
John~C. Baez and Aaron~D. Lauda.
\newblock Higher-dimensional algebra {V}: 2-groups.
\newblock {\em Theory and Applications of Categories}, 12:423--491, 2004,
  math.QA/0307200.

\bibitem{bkp}
R.~Blackwell, G.M. Kelly, and A.J. Power.
\newblock Two-dimensional monad theory.
\newblock {\em Journal of Pure and Applied Algebra}, 59:1--41, 1989.

\bibitem{cheng+lauda}
Eugenia Cheng and Aaron Lauda.
\newblock Higher-dimensional categories: an illustrated guide book, Jun 2004.

\bibitem{fiore}
Thomas~M. Fiore.
\newblock Pseudo limits, biadjoints, and pseudo algebras: Categorical
  foundations of conformal field theory.
\newblock {\em Memoirs of the AMS}, to appear, math.CT/04028298.

\bibitem{h+k}
P.~Hu and I.~Kriz.
\newblock Conformal field theory and elliptic cohomology.
\newblock {\em Advances in Mathematics}, 189(2):325--412, 2004.

\bibitem{j+s}
Andr{\'e} Joyal and Ross Street.
\newblock Braided tensor categories.
\newblock {\em Advances in Mathematics}, 102:20--78, 1993.

\bibitem{kelly}
G.~M. Kelly.
\newblock {\em An Abstract Approach To Coherence}, pages 106--174.
\newblock Number 281 in Lecture Notes in Mathematics. Springer-Verlag, May
  1972.

\bibitem{laplaza}
Miguel~L. Laplaza.
\newblock {\em Coherence for distributivity}, pages 29--65.
\newblock Number 281 in Lecture Notes in Mathematics. Springer-Verlag, May
  1972.

\bibitem{hohc}
Tom Leinster.
\newblock {\em Higher Operads, Higher Categories}.
\newblock Cambridge University Press, 2003, math.CT/0305049.

\bibitem{catwork}
Saunders Mac~Lane.
\newblock {\em Categories for the Working Mathematician}, volume~5 of {\em
  Graduate Texts in Mathematics}.
\newblock Springer-Verlag, second edition, 1998.

\bibitem{modulecat}
Viktor Ostrik.
\newblock Module categories, weak {H}opf algebras and modular invariants.
\newblock {\em Transform. Groups}, 8(2):159--176, 2003, math.QA/0111139.

\bibitem{penon}
J.~Penon.
\newblock Approche polygraphique des $\infty$-cat\'egories non strictes.
\newblock {\em Cahiers de Topologie et G\'eometrie Diff\'erentielle
  Cat\'egoriques}, XL(1):31--80, 1999.

\bibitem{tronin}
S.N. Tronin.
\newblock Abstract clones and operads.
\newblock {\em Siberian Mathematical Journal}, 43(4):746--755, Jul 2002.

\end{thebibliography}
\end{document}